\documentclass{article}

\usepackage{amsmath}
\usepackage{amssymb}

\newtheorem{thm}{Theorem}
\newtheorem{conj}{Conjecture}
\newtheorem{lem}{Lemma}
\newtheorem{prop}{Proposition}
\newtheorem{cor}{Corollary}

\title{Variance of Squarefull Numbers in Short Intervals}
\author{Tsz Ho Chan}
\date{}

\begin{document}
\maketitle

\begin{abstract}
In this paper, we study the variance of the number of squarefull numbers in short intervals. As a result, we are able to prove that, for any $0 < \theta < 1/2$, almost all short intervals $(x, x + x^{1/2 + \theta}]$ contain about $\frac{\zeta(3/2)}{2 \zeta(3)} x^\theta$ squarefull numbers.
\end{abstract}

\section{Introduction and Main Results}

A number $n$ is {\it squarefree} if $p^2 \nmid \; n$ for every $p \mid n$, i.e. all the exponents in its prime factorization are exactly one. In contrast, a number $n$ is {\it squarefull} if $p^2 \mid n$ for every prime $p \mid n$, i.e. all the exponents in its prime factorization are at least two. It is well-known that every squarefull number $n$ can be factored uniquely as $n = a^2 b^3$ for some integer $a$ and squarefree number $b$. Let $Q(x)$ be the number of squarefull numbers which are less than or equal to $x$. 
Then
\begin{equation} \label{squarefull}
Q(x) = \sum_{a^2 b^3 \le x} \mu^2(b)
\end{equation}
where $\mu(n)$ is the M\"{o}bius function. Bateman and Grosswald \cite{BG} proved that
\begin{equation} \label{fullcount}
Q(x) = \frac{\zeta(3/2)}{\zeta(3)} x^{1/2} + \frac{\zeta(2/3)}{\zeta(2)} x^{1/3} + O \bigl( x^{1/6} e^{-C_1 (\log x)^{4/7} (\log \log x)^{-3/7}} \bigr)
\end{equation}
for some absolute constant $C_1 > 0$ where $\zeta(s)$ is the Riemann zeta function. We are interested in squarefull numbers in short intervals, namely
\[
Q(x + y) - Q(x) \; \; \text{ with } \; \; y = o(x).
\]
Bateman and Grosswald's result \eqref{fullcount} gives
\begin{equation} \label{shortcount1}
Q \bigl(x + x^{1/2 + \theta} \bigr) - Q(x) \sim \frac{\zeta(3/2)}{2 \zeta(3)} x^\theta
\end{equation}
for $1/6 = 0.1666... \le \theta < 1/2$. Subsequently, exponential sum, elementary, and differencing methods were used to expand the range of validity for \eqref{shortcount1}. The current best result was obtained by Trifonov \cite{T}, showing that the asymptotic formula \eqref{shortcount1} is true for $19 / 154 = 0.1233... < \theta < 1/2$ (please refer to the same paper for earlier developments and relevant references). As a consequence, for $x^\epsilon \ll H \ll x^{1/2 - \epsilon}$, we have
\begin{equation} \label{shortcount2}
Q \bigl(x + 2 H \sqrt{x} + H^2 \bigr) - Q(x) \ll_\epsilon (1_{H \le x^{19/154}} \cdot x^{19/154} + H) H^\epsilon
\end{equation}
with the indicator function $1_{H \le T} = \left\{ \begin{array}{ll} 1, & \text{ if } H \le T, \\ 0, & \text{ otherwise} \end{array} \right.$. We will use this crude upper bound later.

\bigskip

Inspired by the recent breakthrough work of Gorodetsky, Matom\"{a}ki, Radziwi\l{}\l{} and Rodgers \cite{GMRR} on variance of squarefree numbers in short intervals and arithmetic progressions, we are interested in studying the variance of squarefull numbers in short intervals:
\[
\frac{1}{X} \int_{X}^{2 X} \Big| Q(x + y) - Q(x) - \frac{\zeta(3/2)}{\zeta(3)} (\sqrt{x + y} - \sqrt{x}) \Big|^2 dx.
\]
Our main result is the following.
\begin{thm} \label{thm1}
Given $\epsilon > 0$, $X > 1$ and $X^{3 \epsilon} \le H \le X^{1/2 - \epsilon}$. With $y = 2 \sqrt{x} H + H^2$,
\[
\frac{1}{X} \int_{X}^{2 X} \Big|Q(x + y) - Q(x) - \frac{\zeta(3/2)}{\zeta(3)} H \Big|^2 dx \ll_\epsilon H^{16 / 11 + \epsilon}.
\]
Note that $16 / 11 = 1.454545\ldots$ and $\sqrt{x + y} - \sqrt{x} = H$ (which explains our choice of $y$).
\end{thm}
Define
\[
S := \Bigl\{ X \le x \le 2X :  \Big| Q(x + y) - Q(x) - \frac{\zeta(3/2)}{\zeta(3)} H \Big| >H^{8 / 11 + \epsilon} \Bigr\}.
\]
Then Theorem \ref{thm1} implies $|S| \cdot (H^{8 / 11 + \epsilon})^2 \ll_\epsilon X \cdot H^{16 / 11 + \epsilon}$ where $|S|$ stands for the measure of $S$. Hence, $|S| \ll_\epsilon X / H^\epsilon$ which gives the following.
\begin{cor} \label{cor1}
Given $\epsilon > 0$, $X > 1$ and $X^{3 \epsilon} \le H \le X^{1/2 - \epsilon}$. With $y = 2 \sqrt{x} H + H^2$,
\[
Q(x + y) - Q(x) = \frac{\zeta(3/2)}{\zeta(3)} H + O_\epsilon(H^{8 / 11 + \epsilon}) 
\]
for almost all $x \in [X, 2X]$.
\end{cor}
Since $\epsilon$ can be arbitrarily small, Corollary \ref{cor1} immediately implies
\begin{cor}
For any $0 < \theta < 1/2$, the asymptotic formula \eqref{shortcount1} is true for almost all $x \in [X, 2X]$.
\end{cor}
In both corollaries, almost all means all $x \in [X, 2X]$ apart from a set of measure $o(X)$. Hence, we confirm the expected count \eqref{shortcount1} of squarefull numbers in almost all short intervals. Based on Theorem \ref{thm1}, we suspect the following.
\begin{conj} \label{conj1}
Given $\epsilon > 0$, $X > 1$ and $X^\epsilon \le H \le X^{1/4 - \epsilon}$. With $y = 2 \sqrt{x} H + H^2$,
\[
\frac{1}{X} \int_{X}^{2 X} \Big|Q(x + y) - Q(x) - \frac{\zeta(3/2)}{\zeta(3)} H \Big|^2 dx \sim \frac{4 \zeta(4/3)}{3 \zeta(2)} \int_{0}^{\infty} \Bigl(\frac{\sin \pi y}{\pi y}\Bigr)^2 y^{1/3} dy \cdot H^{2/3}.
\]
\end{conj}
The reason for the upper exponent $1/4 - \epsilon$ is that the secondary main term from \eqref{fullcount} gives
\[
\frac{\zeta(2/3)}{\zeta(2)} (\sqrt[3]{x+y} - \sqrt[3]{x}) \gg \frac{y}{X^{2/3}} \gg \frac{H}{X^{1/6}} \gg H^{1/3}
\]
when $H \ge X^{1/4}$. Unfortunately, we were unable to establish the above conjecture unconditionally for any range of $H$ at the moment. To do that, one needs to improve the error terms in Proposition \ref{prop1} below. However, assuming the Riemann Hypothesis or the Lindel\"{o}f Hypothesis $\zeta(1/2 + it) \ll_\epsilon t^\epsilon$, one has the following conditional result.
\begin{thm} \label{thm2}
Given $\epsilon > 0$ and $X > 1$. Conjecture \ref{conj1} is true for $X^\epsilon \le H \le X^{3 / 14 - \epsilon}$ under the Lindel\"{o}f Hypothesis. Note that $3/14 = 0.2142857 \ldots$.
\end{thm}

One future direction of research is to lower the exponent $16/11$ in Theorem \ref{thm1} to something below $1$. To do that, one may need to include the secondary main term from \eqref{fullcount} and consider
\[
\frac{1}{X} \int_{X}^{2 X} \Big| Q(x + y) - Q(x) - \frac{\zeta(3/2)}{\zeta(3)} (\sqrt{x + y} - \sqrt{x}) - \frac{\zeta(2/3)}{\zeta(2)} (\sqrt[3]{x+y} - \sqrt[3]{x}) \Big|^2 dx
\]
instead. It would also be interesting to see if one can prove the full range for Conjecture \ref{conj1} under the Lindel\"{o}f Hypothesis. Another direction of inquiry is to study the variance of squarefull numbers over arithmetic progressions. However, the method of Gorodetsky, Matom\"{a}ki, Radziwi\l{}\l{} and Rodgers \cite{GMRR} cannot be applied directly as they made use of Poisson summation which does not translate easily into the ``non-linear" situation involving squarefull numbers. It is worth mentioning that Roditty-Gershon \cite{R} recently studied the variance of the number of squarefull polynomials over finite fields in both short intervals and arithmetic progressions settings, and obtained interesting asymptotic results.

\bigskip

Recently, Matom\"{a}ki and Radziwi\l{}\l{} \cite{MR} made an important discovery on multiplicative function $f(n)$ with $-1 \le f(n) \le 1$ by showing that the local average of $f$ on (almost all) short interval $[x, x+h]$ is close to its global average between $[X, 2X]$. It would be interesting to see if similar behavior holds for multiplicative functions supported on squarefull numbers or other sparse sequences.

\bigskip

Our proof of Theorem \ref{thm1} uses different strategies based on various sizes of $H$. For large $H$, one simply applies Bateman and Grosswald's asymptotic formula \eqref{fullcount}. For small $H$, we establish the following two propositions by means of smooth weights as in \cite{GMRR} and some counting arguments.
\begin{prop} \label{prop1}
Given $0 < \epsilon < 0.1$. For $X^{3 \epsilon} \le H \le X^{1/2 - \epsilon}$ and $y = 2 \sqrt{x} H + H^2$, we have
\begin{align*}
\frac{1}{X} \int_{X}^{2X} \Big| Q(x+y) - Q(x) - \frac{\zeta(3/2)}{\zeta(3)} H \Big|^2 dx
=& (1 + O(X^{-\epsilon /2})) 2 H^2  \sum_{b \le X^{1/4}} \frac{\mu^2(b)}{b^3} \sum_{n \ge 1} S\Bigl( \frac{n H}{b^{3/2}} \Bigr)^2 \\
&+ O_\epsilon \Bigl( H^{4/3} + \frac{H^{5/3 + \epsilon}}{X^{1/16}} + \frac{H^{2+\epsilon}}{X^{1/8}} \Bigr).
\end{align*}
Note: We keep the expected main term even though Proposition \ref{prop2} below renders it much smaller than the error term.
\end{prop}
\begin{prop} \label{prop2}
Given $0 < \epsilon < 0.1$. For $1 \le H \le X^{1/2}$, we have
\[
2 H^2 \sum_{b \le X^{1/4}} \frac{\mu^2(b)}{b^3} \sum_{n \ge 1} S\Bigl( \frac{n H}{b^{3/2}} \Bigr)^2= H^{2/3} \cdot \frac{4 \zeta(4/3)}{3 \zeta(2)} \int_{0}^{\infty} S(y)^2 y^{1/3} dy + O_\epsilon \Bigl(\frac{H^{1 + \epsilon/4}}{X^{1/8}} + H^{2/3 - \epsilon/6} \Bigr)
\]
where $S(x) = \frac{\sin \pi x}{\pi x}$. Note: This is truly an asymptotic formula when $H \le X^{3/8 - \epsilon}$.
\end{prop}

This paper is organized as follows. First, we will give explicit construction of certain smooth weights and prove a related lemma on oscillatory integral. We also recall the Kusmin-Landau inequality on exponential sums and a strong form of the Erd\H{o}s-Tur\'{a}n inequality on uniform distribution. Then we will prove Proposition \ref{prop1} by splitting the sum into two major parts. The first part is estimated by means of gaps of fractions while the second part requires a discretization and switching of roles between two variables. This provides another way to approach the problem other than the usage of Dirichlet series and their mean-value theorems. Next, we will prove Proposition \ref{prop2}. Then we will combine everything together to derive Theorem \ref{thm1}. At the end, we will give a sketch of the proof of Theorem \ref{thm2} by highlighting modifications of our arguments through combination with Gorodetsky, Matom\"{a}ki, Radziwi\l{}\l{} and Rodgers' method.

\bigskip

{\bf Notation.} We use $[x]$, $\{x\}$ and $\| x \|$ to denote the integer part of $x$, the fractional part, and the distance between $x$ and the nearest integer respectively. The symbols $f(x) = O(g(x))$, $f(x) \ll g(x)$ and $g(x) \gg f(x)$ are equivalent to $|f(x)| \leq C g(x)$ for some constant $C > 0$. $f(x) \asymp g(x)$ means that $C_1 f(x) \le g(x) \le C_2 f(x)$ for some constants $0 < C_1 < C_2$. $f(x) \sim g(x)$ means that $\lim_{x \rightarrow \infty} f(x)/g(x) = 1$. Finally $f(x) = O_{\lambda} (g(x))$, $f(x) \ll_{\lambda} g(x)$, $g(x) \gg_{\lambda} f(x)$ or $g(x) \asymp_{\lambda} f(x)$ mean that the implicit constant may depend on $\lambda$.
\section{Construction of smooth functions and some lemmas}

Let
\[
u(x) := \left\{ \begin{array}{ll} 1, & \text{ if } x \ge 0, \\ 0, & \text{ if } x < 0, \end{array} \right. \; \; \text{ and } \; \; r_{L}(x) := \left\{ \begin{array}{ll} \frac{1}{L}, & \text{ if } 0 \le x \le L, \\
0, & \text{ otherwise} \end{array} \right.
\]
for some parameter $0 < L \le X$. Recall the definition of the convolution of two integrable real-valued functions $f$ and $g$ as
\[
(f * g)(x) := \int_{-\infty}^{\infty} f(x - u) g(u) du.
\]
We define $u_{k}(x)$ recursively as follows:
\[
u_{0}(x) = u(x), \; \; \text{ and } \; \; u_{k+1}(x) = (u_{k} * r_{L})(x) \text{ for integer } k \ge 0.
\]
One can show that $0 \le u_{k}(x) \le 1$ on $[0, \infty)$, $u_k(x) = 0$ when $x \in (-\infty, 0]$ and $u_{k}(x) = 1$ when $x \in [k L, \infty)$. Also, one can prove inductively that $u_{k}(x)$ is $k-1$ times continuously differentiable for $k \ge 2$ and satisfies the derivative bound $u_{k}^{(i)}(x) \ll_k \frac{1}{L^i}$ for $1 \le i \le k-1$ by using $(f * g)' = f' * g$ repeatedly. In fact, one can derive an explicit formula for $u_{k}(x)$ via Laplace transform $\mathcal{L}$ and its properties as follows:
\begin{align*}
\mathcal{L} \{u_{k}(x)\} =& \mathcal{L}\{u(x)\} \cdot \mathcal{L} \{r_L(x)\}^{k} = \frac{1}{s} \cdot \Bigl[ \frac{1}{L} \cdot \frac{1 - e^{-L s}}{s} \Bigr]^{k} = \frac{1}{L^{k}} \cdot \frac{(1 - e^{-L s})^k}{s^{k+1}}\\
=& \frac{1}{k! L^k} \sum_{i = 0}^{k} (-1)^i \binom{k}{i} \frac{k! e^{-i L s}}{s^{k+1}}, \\
u_{k}(x) =& \mathcal{L}^{-1} \Bigl\{ \frac{1}{k! L^k} \sum_{i = 0}^{k} (-1)^i \binom{k}{i} \frac{k! e^{-i L s}}{s^{k+1}} \Bigr\} = \frac{1}{k! L^k} \sum_{i = 0}^{k} (-1)^i \binom{k}{i} u(x - i L) (x - i L)^{k}.
\end{align*}
For $0 < L < X / (2k)$, we define the smooth weight functions
\[
\sigma_{k,X,L}^-(x) := u_k(x - X) \cdot u_k(2X - x) \; \; \text{ and } \; \; \sigma_{k,X,L}^+(x) := u_k(x - X + kL) \cdot u_k(2X + kL - x).
\]
One can verify that $0 \le \sigma_{k,X,L}^-(x) \le 1$ and $0 \le \sigma_{k,X,L}^+(x) \le 1$ are supported on $[X, 2X]$ and $[X - kL, 2X + k L]$ respectively. Moreover, $\sigma_{k,X,L}^-(x) = 1$ for $x \in [X+kL, 2X-kL]$ and $\sigma_{k,X,L}^+(x) = 1$ for $x \in [X, 2X]$. One also has the derivative conditions:
\[
\frac{d^{i} \sigma_{k, X, L}^{\pm}(x)}{d x^i} \ll_k \frac{1}{L^i} \text{ for } 1 \le i \le k - 1.
\]

We will need the following lemma to show that oscillatory integral of certain smooth function is small.
\begin{lem} \label{lem1}
Suppose $K > 1$ is some large integer. Let $H \ge 1$, $W, Q, L, R > 0$ satisfying $Q R / \sqrt{H} > 1$ and $R L > 1$. Suppose that $w(t)$ is a smooth function with support on $[\alpha, \beta]$ satisfying
\[
w^{(j)}(t) \ll_j W L^{-j} \; \; \text{ for } j = 0, 1, 2, ..., 2K-1.
\]
Suppose $h(t)$ is a smooth function on $[\alpha, \beta]$ such that $|h'(t)| \ge R$ and
\[
h^{(j)}(t) \ll_j H Q^{-j}, \; \; \text{ for } j = 2, 3, ..., 2K- 1.
\]
Then the integral
\[
I = \int_{-\infty}^{\infty} w(t) e^{i h(t)} dt
\]
satisfies
\[
I \ll_K (\beta - \alpha) W [ (Q R / \sqrt{H})^{-K} + (R L)^{-K} ].
\]
\end{lem}

Proof: This is basically Lemma 8.1 in \cite{BKY} but we give the proof here for completeness. Note that we do not require $w$ to be infinitely differentiable or in Schwartz class. Since $w$ has compact support and $h'$ is non-zero, repeated integration by parts gives
\begin{align*}
I =& \int_{-\infty}^{\infty} w(t) d \frac{e^{i h(t)}}{i h'(t)} = \int_{-\infty}^{\infty} \frac{w'(t)}{i h'(t)} e^{i h(t)} dt = \int_{-\infty}^{\infty} \frac{w'(t)}{i h'(t)} d \frac{e^{i h(t)}}{i h'(t)} \\
=& \int_{-\infty}^{\infty} \frac{1}{i h'(t)} \Bigl(\frac{w''(t)}{i h'(t)} - \frac{w'(t) h''(t)}{i h'(t)^2} \Bigr) e^{i h(t)} dt \\
=& \int_{-\infty}^{\infty} \Bigl( -\frac{w''(t)}{h'(t)^2} + \frac{w'(t) h''(t)}{h'(t)^3} \Bigr) d \frac{e^{i h(t)}}{i h'(t)} \\
=& \int_{-\infty}^{\infty} \frac{1}{i h'(t)} \Bigl(- \frac{w'''(t)}{h'(t)^2} + \frac{2 w''(t) h''(t)}{h'(t)^3} + \frac{w''(t) h''(t)}{h'(t)^3} + \frac{w'(t) h'''(t)}{h'(t)^3} - \frac{3 w'(t) h''(t)^2}{h'(t)^4} \Bigr) e^{i h(t)} dt \\
=& \int_{-\infty}^{\infty} \Bigl(- \frac{w'''(t)}{i h'(t)^3} + \frac{3 w''(t) h''(t)}{i h'(t)^4} + \frac{w'(t) h'''(t)}{i h'(t)^4} - \frac{3 w'(t) h''(t)^2}{i h'(t)^5} \Bigr) d \frac{e^{i h(t)}}{i h'(t)} \\
&\cdots \\
=& \int_{-\infty}^{\infty} \sum_{\nu = K}^{2K-1} \sum_{\mu = 1}^{\nu} \frac{w^{(\mu)}(t)}{h'(t)^{\nu}} \sum_{2 \gamma_2 + 3 \gamma_3 + \cdots + \nu \gamma_{\nu} = \nu - \mu} c_{\nu, \mu, \gamma_2, \gamma_3, \cdots, \gamma_{\nu}} h^{(2)}(t)^{\gamma_2} h^{(3)}(t)^{\gamma_3} \cdots h^{(\nu)}(t)^{\gamma_{\nu}} e^{i h(t)} dt
\end{align*}
for some absolute coefficients $c_{\nu, \mu, \gamma_2, \gamma_3, \cdots, \gamma_{\nu}} \in \mathbb{C}$. Hence,
\[
|I| \ll_K (\beta - \alpha) W \sum_{\nu = K}^{2K-1} R^{-\nu} \sum_{\mu = 1}^{\nu} L^{-\mu} \frac{H^{(\nu - \mu) / 2}}{Q^{\nu - \mu}}
\]
which gives the lemma.

\bigskip

We also need two standard results on exponential sums and uniform distribution of sequences.

\begin{lem}[Kusmin-Landau inequality] \label{kusmin-landau}
Let $f : (a, b] \rightarrow \mathbb{R}$ be a continuously differentiable function and, for some real number $0 < \lambda < 1$, $\| f'(x) \| \ge \lambda$ for all $x \in (a, b]$. Then
\[
\sum_{a < n \le b} e ( f(n) ) \ll \frac{1}{\lambda}
\]
where the implicit constant does not depend on $f$. Here $e(u) = e^{2 \pi i u}$.
\end{lem}

Proof: See Theorem 2.1 in \cite{GK} for example.

\begin{lem}[Strong form of Erd\H{o}s-Tur\'{a}n inequality] \label{erdos-turan}
Let $x_1, x_2, x_3, \ldots$ be a sequence of real numbers and let $0 \le \alpha < \beta < 1$. Define the discrepancy
\[
D_N(\alpha, \beta) := \# \{ 1 \le n \le N : \alpha \le \{ x_n \} \le \beta \} - N (\beta - \alpha).
\]
Then
\[
D_N(\alpha, \beta) \le \frac{N}{K+1} + 2 \sum_{k = 1}^{K} \Bigl( \frac{1}{K+1} + \min \Bigl( \beta - \alpha, \frac{1}{\pi k} \Bigr) \Bigr) \bigg| \sum_{n = 1}^{N} e(k u_n) \bigg| 
\]
for any positive integers $N$ and $K$.
\end{lem}

Proof: See Theorem 1 in Chapter 1 of \cite{M} for example.


\section{Proof of Proposition \ref{prop1}: Initial Manipulations}

By \eqref{squarefull}, we can rewrite
\begin{equation} \label{z1}
I := \frac{1}{X} \int_{X}^{2X} \Big| Q(x+y) - Q(x) - \frac{\zeta(3/2)}{\zeta(3)} H \Big|^2 dx = \frac{1}{X} \int_{X}^{2X} \Big| \sum_{x < a^2 b^3 \le x + y} \mu^2(b) - \frac{\zeta(3/2)}{\zeta(3)} H \Big|^2 dx 
\end{equation}
where $y = 2 \sqrt{x} H + H^2$ with $H \le X^{1/2 - \epsilon}$. Observe that
\begin{equation} \label{zeta}
\sum_{b = 1}^{\infty} \frac{\mu^2(b)}{b^s} = \prod_{p} \Bigl(1 + \frac{1}{p^s} \Bigr) = \prod_{p} \frac{(1 - \frac{1}{p^{2s}})}{(1 - \frac{1}{p^s})} = \frac{\zeta(s)}{\zeta(2s)}
\end{equation}
for $\Re{s} > 1$. To study $I$, we separate those squarefull numbers with $b \le X^{1/4}$ and $b > X^{1/4}$, and further break up the latter case into dyadic ranges $a \in [A_i, 2 A_i)$ where $A_i = 2^{i-1}$ with $1 \le i \le \frac{\log 2X}{8 \log 2}$. By Cauchy-Schwarz inequality, we have
\begin{align} \label{Ientire}
I =& \frac{1}{X} \int_{X}^{2X} \bigg| \mathop{\sum_{x < a^2 b^3 \le x + y}}_{b \le X^{1/4}} \mu^2(b) + \sum_{i} \mathop{\sum_{x < a^2 b^3 \le x + y}}_{b > X^{1/4}, \, A_i \le a < 2 A_i} \mu^2(b)  - H \sum_{b = 1}^{\infty} \frac{\mu^2(b)}{b^{3/2}} \bigg|^2 dx \nonumber \\
=& \frac{1}{X} \int_{X}^{2X} \bigg| \mathop{\sum_{x < a^2 b^3 \le x + y}}_{b \le X^{1/4}} \mu^2(b) - H \sum_{b \le X^{1/4}} \frac{\mu^2(b)}{b^{3/2}} + \sum_{i} \mathop{\sum_{x < a^2 b^3 \le x + y}}_{b > X^{1/4}, \, A_i \le a < 2 A_i} \mu^2(b)  + O \Bigl( \frac{H}{X^{1/8}} \Bigr) \bigg|^2 dx \nonumber \\
=& J_1 + J_2 + O \Bigl( J_1^{1/2} J_2^{1/2} + \frac{J_1^{1/2} H}{X^{1/8}} + \frac{J_2^{1/2} H}{X^{1/8}} + \frac{H^2}{X^{1/4}} \Bigr)
\end{align}
where
\begin{equation} \label{J1}
J_1 := \frac{1}{X} \int_{X}^{2X} \bigg| \mathop{\sum_{x < a^2 b^3 \le x + y}}_{b \le X^{1/4}} \mu^2(b) - H \sum_{b \le X^{1/4}} \frac{\mu^2(b)}{b^{3/2}} \biggr|^2 dx,
\end{equation}
and
\begin{equation} \label{J2}
J_2 := \frac{1}{X} \int_{X}^{2X} \bigg| \sum_{i} \mathop{\sum_{x < a^2 b^3 \le x + y}}_{b > X^{1/4}, \, A_i \le a < 2 A_i} \mu^2(b) \biggr|^2 dx.
\end{equation}

\section{Proof of Proposition \ref{prop1}: Estimation of $J_1$}

We are going to study a smoothed version of $J_1$ as defined in \eqref{J1}. Let $L = X^{1 - \epsilon/2}$ and $K$ be some sufficiently large integer (in terms of $\epsilon$). Consider
\begin{equation} \label{z-smallsmooth}
\frac{1}{X} \int_{-\infty}^{\infty} \sigma_{K, X, L}^{\pm} (x) \bigg| \mathop{\sum_{x < a^2 b^3 \le x + y}}_{b \le X^{1/4}} \mu^2(b) - H \sum_{b \le X^{1/4}} \frac{\mu^2(b)}{b^{3/2}} \bigg|^2 dx.
\end{equation}
We have
\begin{align} \label{ab-sum}
\mathop{\sum_{x < a^2 b^3 \le x + y}}_{b \le X^{1/4}} \mu^2(b) =& \sum_{b \le X^{1/4}} \mu^2(b) \sum_{\sqrt{x / b^3} < a \le \sqrt{(x + y) / b^3}} 1 \nonumber \\
=& \sum_{b \le X^{1/4}} \mu^2(b) \Bigl[ \sqrt{\frac{x + y}{b^3}} - \sqrt{\frac{x}{b^3}} - \psi \Bigl(\sqrt{\frac{x + y}{b^3}} \Bigr) + \psi \Bigl(\sqrt{\frac{x}{b^3}} \Bigr) \Bigr] \nonumber \\
=& H \sum_{b \le X^{1/4}} \frac{\mu^2(b)}{b^{3/2}} - \sum_{b \le X^{1/4}} \mu^2(b) \Bigl[ \psi \Bigl(\sqrt{\frac{x + y}{b^3}} \Bigr) - \psi \Bigl(\sqrt{\frac{x}{b^3}} \Bigr) \Bigr] 
\end{align}
where $\psi(u) = u - [u] - 1/2$. For $\psi(u)$, we apply the truncated Fourier expansion
\begin{equation} \label{frac-fourier}
\psi(u) = - \frac{1}{2 \pi i} \sum_{0 < |n| \le N} \frac{1}{n} e(n u) + O \Bigl( \min \bigl(1, \frac{1}{N \| u \|} \bigr) \Bigr).
\end{equation}
We take $N = X^{10}$ and put \eqref{frac-fourier} into \eqref{ab-sum}. The arising error term is $O(1 / X^4)$ unless $\| \sqrt{(x + y) / b^3} \| < X^{-5}$ or $\| \sqrt{x / b^3} \| < X^{-5}$. When $\| \sqrt{(x + y) / b^3} \| < X^{-5}$ or $\| \sqrt{x / b^3} \| < X^{-5}$, we have $|\sqrt{(x + y) / b^3} - a_1| < X^{-5}$ or $| \sqrt{x / b^3} - a_2| < X^{-5}$ for some integers $a_1, a_2$ which implies
\begin{equation} \label{support}
\Big| \frac{x+y}{b^3} - a_1^2 \Bigr|, \; \Big| \frac{x}{b^3} - a_2^2 \Bigr| \ll \frac{1}{X^{4.5}} \; \; \text{ or } \; \; |x+y - a_1^2 b^3|, \; |x - a_2^2 b^3| \ll \frac{1}{X^{3.5}}
\end{equation}
as $b \le X^{1/4}$. Let $E(x)$ be the error in the second sum in \eqref{ab-sum} coming from the error term of \eqref{frac-fourier} when $\| \sqrt{(x + y) / b^3} \| < X^{-5}$ or $\| \sqrt{x / b^3} \| < X^{-5}$ for some $b \le X^{1/4}$. As there are $O(\sqrt{X})$ squarefull numbers in the interval $[X, 3X]$ by \eqref{fullcount},
\begin{equation} \label{E}
E(x) \ll X^{1/4}, \; \; \text{ and } \; \; E(x) \text{ is supported on } O(\sqrt{X}) \text{ intervals of length } O\Bigl(\frac{1}{X^{3.5}} \Bigr) \text{ each}
\end{equation}
by \eqref{support}. Thus, we can replace \eqref{z-smallsmooth} by
\begin{align} \label{z-smallsmooth2}
\frac{1}{4 \pi^2 X} & \int_{-\infty}^{\infty} \sigma_{K, X, L}^{\pm} (x) \bigg| \sum_{b \le X^{1/4}} \mu^2(b) \sum_{0 < |n| \le N} \frac{1}{n} e \Bigl(\frac{n \sqrt{x}}{b^{3/2}} \Bigl) \Bigl[1 - e \Bigl( \frac{n H}{b^{3/2}} \Bigr) \Bigr] + E(x) + O\Bigl(\frac{1}{X^4} \Bigr) \bigg|^2 dx \nonumber \\
=& \frac{1}{X} \Bigl[J_{\pm} + O \Bigl(\frac{J_{\pm}^{1/2}}{X^{5/4}} + \frac{1}{X^{5/2}}\Bigr) \Bigr]
\end{align}
by Cauchy-Schwarz inequality and $\int_{-\infty}^{\infty} \sigma_{K, X, L}^{\pm} (x) |E(x) + O(1/X^4)|^2 dx \ll X^{- 5/2}$ via \eqref{E}. Here
\begin{equation} \label{z-smallsmooth2.5}
J_{\pm} := \frac{1}{4 \pi^2} \int_{-\infty}^{\infty} \sigma_{K, X, L}^{\pm} (x) \bigg| \sum_{b \le X^{1/4}} \mu^2(b) \sum_{0 < |n| \le N} \frac{1}{n} e \Bigl(\frac{n \sqrt{x}}{b^{3/2}} \Bigl) \Bigl[1 - e \Bigl( \frac{n H}{b^{3/2}} \Bigr) \Bigr] \bigg|^2 dx.
\end{equation}
Note that, since the characteristic function of the interval $[X, 2X]$ is between $\sigma_{K,X,L}^-$ and $\sigma_{K,X,L}^+$, we have
\begin{equation} \label{newintermediate}
\frac{1}{X} \Bigl[J_{-} + O \Bigl(\frac{J_{-}^{1/2}}{X^{5/4}} + \frac{1}{X^{5/2}}\Bigr) \Bigr] \le J_1 \le \frac{1}{X} \Bigl[J_{+} + O \Bigl(\frac{J_{+}^{1/2}}{X^{5/4}} + \frac{1}{X^{5/2}}\Bigr) \Bigr].
\end{equation}
Expanding \eqref{z-smallsmooth2.5} out,
\begin{align} \label{z-smallsmooth3}
J_{\pm} =\frac{1}{4 \pi^2} & \sum_{b_1, b_2 \le X^{1/4}} \sum_{0 < |n_1|, |n_2| \le N} \mu^2(b_1) \mu^2(b_2) \frac{1}{n_1 n_2} \nonumber \\
&\times \Bigl[1 - e \Bigl( \frac{n_1 H}{b_1^{3/2}} \Bigr) \Bigr] \overline{\Bigl[1 - e \Bigl( \frac{n_2 H}{b_2^{3/2}} \Bigr) \Bigr]} \int_{-\infty}^{\infty} \sigma_{K, X, L}^{\pm} (x) e \Bigl(x^{1/2} \Bigl(\frac{n_1}{b_1^{3/2}} - \frac{n_2}{b_2^{3/2}} \Bigr) \Bigr) dx.
\end{align}
By Lemma \ref{lem1} with $H = X^\epsilon$, $W = 1$, $Q = X$, $L = X^{1 - \epsilon/2}$ and $R = 1/X^{1-\epsilon}$, the contribution from those terms with $| \frac{n_1}{b_1^{3/2}} - \frac{n_2}{b_2^{3/2}} | > \frac{1}{X^{1/2 - \epsilon}}$ is
\begin{equation} \label{off}
\ll X^{- K \epsilon / 2} \ll \frac{1}{X^5}
\end{equation}
by picking $K = [20 / \epsilon]$ for example. Hence, we may restrict the sum in \eqref{z-smallsmooth3} to those integers for which
\begin{equation*}
\bigg| \frac{n_1}{b_1^{3/2}} - \frac{n_2}{b_2^{3/2}} \bigg| \le \frac{1}{X^{1/2 - \epsilon}}.
\end{equation*}
We separate those $(n_1, n_2, b_1, b_2)$ for which $\frac{n_1}{b_1^{3/2}} = \frac{n_2}{b_2^{3/2}}$ (i.e. $n_1^2 b_2^3 = n_2^2 b_1^3$) and those for which $\frac{n_1}{b_1^{3/2}} \neq \frac{n_2}{b_2^{3/2}}$. In the first case, we must have $b_1 = b_2 = b$ and $n_1 = n_2 = n$ which contributes
\begin{align*}
\frac{\overline{\sigma_{K,X,L}^{\pm}}}{4 \pi^2} \sum_{b \le X^{1/4}} \mu^2(b) \sum_{n \neq 0} \frac{1}{n^2}
\Big|1 - e \Bigl( \frac{n H}{b^{3/2}} \Bigr) \Big|^2 + O \Bigl(\frac{1}{X^5}\Bigr)
\end{align*}
where the error term comes from adding $|n| > N = X^{10}$ and
\begin{equation} \label{mean-sigma}
\overline{\sigma_{K,X,L}^{\pm}} = \int_{-\infty}^{\infty} \sigma_{K,X,L}^{\pm}(x) dx = X + O_\epsilon(X^{1- \epsilon/2}).
\end{equation}
As
\[
\Big|1 - e \Bigl( \frac{n H}{b^{3/2}} \Bigr) \Big| = 2 \Big| \sin \Bigl( \frac{n \pi H}{b^{3/2}} \Bigr) \Big|,
\]
the diagonal terms ($\frac{n_1}{b_1^{3/2}} = \frac{n_2}{b_2^{3/2}}$) from \eqref{z-smallsmooth3} contribute 
\begin{equation} \label{diagonal}
(1 + O(X^{-\epsilon /2})) X H^2  \sum_{b \le X^{1/4}} \frac{\mu^2(b)}{b^3} \sum_{n \neq 0} S\Bigl( \frac{n H}{b^{3/2}} \Bigr)^2 + O \Bigl(\frac{1}{X^5}\Bigr)
\end{equation}
where $S(x) = \frac{\sin \pi x}{\pi x}$. It remains to show that the contribution from those non-diagonal terms with $n_1^2 b_2^3 \neq n_2^2 b_1^3$ is also acceptable. Splitting $n_i$ and $b_i$ dyadically, we need to bound
\begin{equation} \label{bound2}
\min \Bigl(\frac{1}{N_1}, \frac{H}{B_1^{3/2}} \Bigr) \min \Bigl(\frac{1}{N_2}, \frac{H}{B_2^{3/2}} \Bigr) \# \Bigl\{ (b_1, b_2, n_1, n_2) : b_i \sim B_i, n_i \sim N_i, 0 < \Big| \frac{n_1}{b_1^{3/2}} - \frac{n_2}{b_2^{3/2}} \Big| \le \frac{1}{X^{1/2 - \epsilon}} \Bigr\}
\end{equation}
for any $B_1, B_2 \le X^{1/4}$ and $N_1, N_2 \le N$. We have $\frac{B_1^{3/2}}{X^{1/2 - \epsilon}}, \frac{B_2^{3/2}}{X^{1/2 - \epsilon}} \ll \frac{1}{X^{1/8 - \epsilon}} < 0.1$. So, we can assume that $N_1 \ge \frac{10 B_1^{3/2}}{X^{1/2 - \epsilon}}$ and $N_2 \ge \frac{10 B_2^{3/2}}{X^{1/2 - \epsilon}}$. Thus,
\begin{equation} \label{dioph-ineq}
\Big| \frac{n_1}{b_1^{3/2}} - \frac{n_2}{b_2^{3/2}} \Big| \le \frac{1}{X^{1/2 - \epsilon}} \; \; \text{ or } \; \; |n_1^2 b_2^3 - n_2^2 b_1^3| \ll \frac{N_1 B_1^{3/2} B_2^{3}}{X^{1/2 - \epsilon}} \asymp \frac{N_1^{1/2} N_2^{1/2} B_1^{9/4} B_2^{9/4}}{X^{1/2 - \epsilon}}
\end{equation}
has no solution unless $N_1 B_2^{3/2} \asymp N_2 B_1^{3/2}$. The following counting argument was inspired by Heath-Brown's work in \cite{H} but we need to take weights into account. Let $n = \gcd(n_1, n_2)$ and $b = \gcd(b_1, b_2)$. We can rewrite
\[
n_1 = n n_1', \; \; n_2 = n n_2', \; \; b_1 = b b_1', \; \; b_2 = b b_2'.
\]
Inequality \eqref{dioph-ineq} implies
\begin{equation*}
\Big| \frac{n_1'}{n_2'} - \frac{b_1'^{3/2}}{b_2'^{3/2}} \Big| \ll \frac{B_1^{3/2}}{N_2 X^{1/2 - \epsilon}}.
\end{equation*}
Note that $|\frac{n_1'}{n_2'} - \frac{n_1''}{n_2''}| \gg \frac{1}{(N_2 / n)^2}$ for distinct fractions. Hence, for fixed greatest common divisors $b$ and $n$, and fixed integers $b_1'$ and $b_2'$, the number of $n_1'$ and $n_2'$ satisfying \eqref{dioph-ineq} is
\[
\ll \frac{B_1^{3/2} / (N_2 X^{1/2 - \epsilon})}{1 / (N_2 / n)^2} + 1 = \frac{B_1^{3/2} N_2}{n^2 X^{1/2 - \epsilon}} + 1.
\]
Thus, the quantity in \eqref{bound2} is
\begin{align} \label{bound3}
&\ll \sum_{b \le 2 B_1} \sum_{n \le 2 N_2} \sum_{b_1' \ll B_1 / b} \; \sum_{b_2' \ll B_2 / b} \min \Bigl(\frac{1}{N_1}, \frac{H}{B_1^{3/2}} \Bigr) \min \Bigl(\frac{1}{N_2}, \frac{H}{B_2^{3/2}} \Bigr) \Bigl( \frac{B_1^{3/2} N_2}{n^2 X^{1/2 - \epsilon}} + 1 \Bigr) \nonumber \\
&\ll \sum_{b} \sum_{n} \sum_{b_1'} \sum_{b_2' } \frac{H}{B_1^{3/2}} \frac{1}{N_2} \frac{B_1^{3/2} N_2}{n^2 X^{1/2 - \epsilon}} + \sum_{b} \sum_{n} \sum_{b_1' \ll B_1 / b} \sum_{b_2' \ll B_2 / b} \frac{1}{N_1^{1/3}} \frac{H^{2/3}}{B_1} \frac{1}{N_2^{1/3}} \frac{H^{2/3}}{B_2} \nonumber \\
&\ll \frac{H B_1 B_2}{X^{1/2 - \epsilon}} + H^{4/3} \ll H X^\epsilon + H^{4/3}
\end{align}
as long as $n \le (N_1 N_2)^{1/3}$. It remains to deal with the situation when $n > (N_1 N_2)^{1/3}$. From \eqref{dioph-ineq}, we also have
\[
\Big| \frac{n_1'^{2/3}}{n_2'^{2/3}} - \frac{b_1'}{b_2'} \Big| \ll \frac{B_1 B_2^{1/2}}{N_2 X^{1/2 - \epsilon}}.
\]
Note that $| \frac{b_1'}{b_2'} - \frac{b_1''}{b_2''} | \gg \frac{1}{(B_2 / b)^2}$. Hence, for fixed greatest common divisors $b$ and $n$, and fixed integers $n_1'$ and $n_2'$, the number of $b_1'$ and $b_2'$ satisfying \eqref{dioph-ineq} is
\[
\ll \frac{ B_1 B_2^{1/2} / (N_2 X^{1/2 - \epsilon})}{1 / (B_2 / b)^2} + 1 = \frac{ B_1 B_2^{5/2}}{b^2 N_2 X^{1/2 - \epsilon}} + 1.
\]
We claim that $b \le (B_1 B_2)^{1/2}$. For otherwise, the inequality $n^2 b^3 \le |n_1^2 b_2^3 - n_2^2 b_1^3|$ together with \eqref{dioph-ineq} would imply
\[
(N_1 N_2)^{2/3} (B_1 B_2)^{3/2} < n^2 b^3 \ll \frac{N_1^{1/2} N_2^{1/2} B_1^{9/4} B_2^{9/4}}{X^{1/2 - \epsilon}} \; \; \text{ or } \; \; X^{1/2 - \epsilon} \ll (B_1 B_2)^{3/4}.
\]
This is impossible as $B_1, B_2 \le X^{1/4}$. Also, note that $n_1' n_2' \le \frac{N_1 N_2}{n^2} \le \frac{(N_1 N_2)^{2/3}}{n}$ as $n > (N_1 N_2)^{1/3}$. Therefore, the quantity in \eqref{bound2} is
\begin{align} \label{bound4}
&\ll \sum_{b \le (B_1 B_2)^{1/2}} \sum_{n \le 2 N_2} \sum_{n_1' \ll N_1 / n} \; \sum_{n_2' \ll N_2 / n} \min \Bigl(\frac{1}{N_1}, \frac{H}{B_1^{3/2}} \Bigr) \min \Bigl(\frac{1}{N_2}, \frac{H}{B_2^{3/2}} \Bigr) \Bigl( \frac{ B_1 B_2^{5/2}}{b^2 N_2 X^{1/2 - \epsilon}} + 1  \Bigr) \nonumber \\
&\ll \sum_{b} \sum_{n} \sum_{n_1'} \sum_{n_2' } \frac{1}{N_1} \frac{H}{B_2^{3/2}} \frac{B_1 B_2^{5/2}}{b^2 N_2 X^{1/2 - \epsilon}} + \sum_{b \le (B_1 B_2)^{1/2}} \sum_{n} \sum_{n_1' n_2' \le (N_1 N_2)^{2/3} / n} \frac{1}{N_1^{2/3}} \frac{H^{1/3}}{B_1^{1/2}} \frac{1}{N_2^{2/3}} \frac{H^{1/3}}{B_2^{1/2}} \nonumber \\
&\ll \frac{H B_1 B_2}{X^{1/2 - \epsilon}} + H^{2/3} X^\epsilon \ll H X^\epsilon
\end{align}
since $B_1, B_2 \le X^{1/4}$. Combining \eqref{z-smallsmooth2.5}, \eqref{newintermediate}, \eqref{z-smallsmooth3}, \eqref{off}, \eqref{mean-sigma}, \eqref{diagonal}, \eqref{bound2}, \eqref{bound3} and \eqref{bound4}, we have
\begin{equation} \label{Jfinal}
J_1 = (1 + O(X^{-\epsilon /2})) 2 H^2  \sum_{b \le X^{1/4}} \frac{\mu^2(b)}{b^3} \sum_{n \ge 1} S\Bigl( \frac{n H}{b^{3/2}} \Bigr)^2 + O_\epsilon ( H^{4/3} )
\end{equation}
when $H \ge X^{3 \epsilon}$. 

\section{Proof of Proposition \ref{prop1}: Estimation of $J_2$}

Recall $A_i = 2^{i-1}$ and $1 \le i \le \frac{\log 2X}{8 \log 2}$. Replacing $\mu(b)^2$ by $1$ and dropping the condition on $b$ in \eqref{J2}, we have
\[
J_2 \ll \frac{1}{X} \int_{X}^{2X} \biggl| \sum_{i} \mathop{\sum_{x < a^2 b^3 \le x + y}}_{A_i \le a < 2 A_i} 1 \biggr|^2 dx.
\]
We will perform a discretization. For any integer $m$ between $X - 1$ and $2X$ and any real number $x \in [m, m+1)$, observe that
\[
\sum_{i} \mathop{\sum_{x < a^2 b^3 \le x + 2 H \sqrt{x} + H^2}}_{A_i \le a < 2 A_i} 1 = \sum_{i} \mathop{\sum_{m < a^2 b^3 \le m + 2 H \sqrt{m} + H^2}}_{A_i \le a < 2 A_i} 1
\]
unless there is a squarefull number $s$ in the interval $( (\sqrt{m} + H)^2, (\sqrt{m+1} + H)^2 ]$ (called such $m$ {\it bad}). Such bad $m$ satisfies $(\sqrt{s} - H)^2 - 1 < m < (\sqrt{s} - H)^2$ for some squarefull number $s$ and
\[
0 \le \sum_{i} \mathop{\sum_{x < a^2 b^3 \le x + 2 H \sqrt{x} + H^2}}_{A_i \le a < 2 A_i} 1 - \sum_{i} \mathop{\sum_{m < a^2 b^3 \le m + 2 H \sqrt{m} + H^2}}_{A_i \le a < 2 A_i} 1 \le 3
\]
as the interval $( (\sqrt{m} + H)^2, (\sqrt{m+1} + H)^2 ]$ has length less than $3$. Hence,
\begin{align} \label{J2new}
J_2 \ll& \frac{1}{X} \mathop{\sum_{X-1 \le m \le 2X}}_{m \text{ not bad}} \biggl| \sum_{i} \mathop{\sum_{m < a^2 b^3 \le m + 2H \sqrt{m} + H^2}}_{A_i \le a < 2 A_i} 1 \biggr|^2 \nonumber \\
&+ \frac{1}{X} \mathop{\sum_{X-1 \le m \le 2X}}_{m \text{ bad}} \biggl|3 + \sum_{i} \mathop{\sum_{m < a^2 b^3 \le m + 2H \sqrt{m} + H^2}}_{A_i \le a < 2 A_i} 1 \biggr|^2 \nonumber \\
\ll_\epsilon& \frac{1}{X}  \biggl( \sum_{X-1 \le m \le 2X} \biggl| \sum_{i} \mathop{\sum_{m < a^2 b^3 \le m + 2H \sqrt{m} + H^2}}_{A_i \le a < 2 A_i} 1 \biggr|^2 + \sqrt{X} (1_{H \le X^{19/154}} \cdot X^{19/154} + H) H^\epsilon \biggr)
\end{align}
as there are $O(\sqrt{X})$ squarefull numbers in $[X, 4X]$ (i.e. at most $O(\sqrt{X})$ bad $m$'s) together with \eqref{shortcount2}. So, it remains to study the second moment of the counting function
\[
C_{m} := \sum_{i} \mathop{\sum_{m < a^2 b^3 \le m + 2H \sqrt{m} + H^2}}_{A_i \le a < 2 A_i} 1.
\]
For fixed $a$, observe that  $m < a^2 b^3 \le m + 2H \sqrt{m} + H^2$ is equivalent to
\[
\frac{m^{1/3}}{a^{2/3}} < b \le \frac{m^{1/3}}{a^{2/3}} + \Bigl(\frac{(\sqrt{m} + H)^{2/3} - \sqrt{m}^{2/3}}{a^{2/3}} \Bigr), 
\]
or $m$ is in the set
\[
\mathcal{S}_a := \Bigl\{ X-1 \le m \le 2X : 1 - \Bigl(\frac{(\sqrt{m} + H)^{2/3} - \sqrt{m}^{2/3}}{a^{2/3}} \Bigr) \le \Bigl\{ \frac{m^{1/3}}{a^{2/3}} \Bigr\} < 1 \Bigr\}.
\]
By Lemma \ref{erdos-turan}, we have
\begin{align} \label{ET}
\sum_{X-1 \le m \le 2X} & \Bigl[ C_{m} - \sum_{i} \sum_{A_i \le a < 2 A_i} \frac{(\sqrt{m} + H)^{2/3} - \sqrt{m}^{2/3}}{a^{2/3}} \Bigr] \nonumber \\ 
&= \sum_{i} \sum_{A_i \le a < 2 A_i} \biggl[ \mathop{\sum_{X - 1 \le m \le 2X}}_{m \in \mathcal{S}_a} 1 - \sum_{X-1 \le m \le 2X}  \frac{(\sqrt{m} + H)^{2/3} - \sqrt{m}^{2/3}}{a^{2/3}} \biggr] \nonumber \\
&\ll \sum_{i} \sum_{A_i \le a < 2 A_i} \biggl[ \frac{X}{K_i} + \sum_{k = 1}^{K_i} \Bigl( \frac{1}{K_i} + \min \Bigl( \frac{H}{a^{2/3} X^{1/6}}, \frac{1}{k} \Bigr) \Bigr) \Big| \sum_{X - 1 \le m \le 2X} e \Bigl( \frac{k m^{1/3}}{a^{2/3}} \Bigr) \Big| \biggr]
\end{align}
as $\frac{(\sqrt{m} + H)^{2/3} - \sqrt{m}^{2/3}}{a^{2/3}} \asymp \frac{H}{a^{2/3} X^{1/6}}$.
We pick $K_i = [A_i^{2/3} (X - 1)^{2/3}]$. By Lemma \ref{kusmin-landau}, the innermost exponential sum is $\ll a^{2/3} X^{2/3} / k$. By splitting the sum over $k$ into $1 \le k \le a^{2/3} X^{1/6} / H$ and $a^{2/3} X^{1/6} / H < k \le K_i$, one can show that the right hand side of \eqref{ET} is
\[
\ll \sum_{i} ( A_i^{1/3} X^{1/3} + H A_i X^{1/2} \log X + A_i^{5/3} X^{2/3} ) \ll H X^{5/8} \log X + X^{7/8}
\]
and, hence,
\begin{equation} \label{almostC}
\sum_{X - 1 \le m \le 2X} C_{m} \ll \sum_{i} H A_i^{1/3} X^{5/6} + H X^{5/8}  \log X + X^{7/8} \ll H X^{7/8}
\end{equation}
as $A_i \ll X^{1/8}$. By \eqref{shortcount2}, we have $C_{m} \ll_\epsilon (1_{H \le X^{19/154}} \cdot X^{19 / 154} + H) H^\epsilon$. Putting this into \eqref{almostC} and \eqref{J2new}, we obtain
\begin{align*}
\sum_{X - 1 \le m \le 2X} C_{m}^2 &\ll_\epsilon (1_{H \le X^{19/154}} \cdot X^{19 / 154} + H) H^\epsilon \sum_{X - 1 \le m \le 2X} C_{m} \\
&\ll (1_{H \le X^{19/154}} \cdot X^{19 / 154} + H) H^{1 + \epsilon} X^{7/8}
\end{align*}
and
\begin{equation} \label{finalC}
J_2 \ll_\epsilon (1_{H \le X^{19/154}} \cdot X^{19 / 154} + H) \frac{H^{1 + \epsilon}}{X^{1/8}} \ll 1_{H \le X^{19/154}} \cdot \frac{H^{1 + \epsilon}}{X^{1/616}} + \frac{H^{2 + \epsilon}}{X^{1/8}}
\end{equation}
when $X^{3 \epsilon} \le H \le X^{1/2 - \epsilon}$.

\section{Finishing the proof of Proposition \ref{prop1}}

From \eqref{Jfinal} and Proposition \ref{prop2}, we have $J_1 \ll_\epsilon H^{4/3}$. From \eqref{finalC}, one can check that
\[
J_2 \ll_\epsilon \left\{ \begin{array}{ll} H^{1 + \epsilon} / X^{1/616}, & \text{ if } X^{3 \epsilon} \le H \le X^{19/154}, \\ 
H^{2+\epsilon} / X^{1/8}, & \text{ if }  X^{19/154} < H \le X^{1/2 - \epsilon} \end{array} \right.
\]
Putting these estimates of $J_1$ and $J_2$ into \eqref{Ientire}, we have
\begin{align*}
\frac{1}{X} \int_{X}^{2X} \Big| Q(x+y) - Q(x) - \frac{\zeta(3/2)}{\zeta(3)} H \Big|^2 dx
=& (1 + O(X^{-\epsilon /2})) 2 H^2  \sum_{b \le X^{1/4}} \frac{\mu^2(b)}{b^3} \sum_{n \ge 1} S\Bigl( \frac{n H}{b^{3/2}} \Bigr)^2 \\
&+ O_\epsilon \Bigl( H^{4/3} + \frac{H^{5/3 + \epsilon}}{X^{1/16}} + \frac{H^{2+\epsilon}}{X^{1/8}} \Bigr)
\end{align*}
which gives Proposition \ref{prop1}.

\section{Proof of Proposition \ref{prop2}}

In this section, we shall focus on finding an asymptotic formula for the diagonal contribution:
\begin{equation} \label{diagonal-main}
2 H^2 \sum_{b \le X^{1/4}} \frac{\mu^2(b)}{b^3} \sum_{n \ge 1} S\Bigl( \frac{n H}{b^{3/2}} \Bigr)^2 \; \; \text{ where } \; \; S(x) = \frac{\sin \pi x}{\pi x}.
\end{equation}
As in \cite{GMRR}, we consider a smoother sum first. Fix $\epsilon > 0$ and $K_0 > 0$. Suppose that $w: \mathbb{R} \rightarrow \mathbb{C}$ satisfies
\begin{equation} \label{deriv}
|w^{(k)}(y)| \le K_0 \frac{H^{l \epsilon / 4}}{(1 + |y|)^l}
\end{equation}
for all $k, l \in \{0, 1, 2, 3, 4 \}$. We claim that
\begin{equation} \label{smoothW}
2 H^2 \sum_{b \le X^{1/4}} \frac{\mu^2(b)}{b^3} \sum_{n \ge 1} \Big| w\Bigl( \frac{n H}{b^{3/2}} \Bigr) \Big|^2 = H^{2/3} \frac{4 \zeta(4/3)}{3 \zeta(2)} \int_{0}^{\infty} |w(y)|^2 y^{1/3} dy + O_\epsilon\Bigl(\frac{H^{1 + \epsilon/4}}{X^{1/8}} + H^{1/3+ 2\epsilon}\Bigr).
\end{equation}

Proof of claim: This is very similar to Lemma 8 in \cite{GMRR}. Firstly, we complete the sum over all $b$. By \eqref{deriv} for $k = 0$ and $l = 0, 1$, we have
\[
\sum_{n \ge 1} |w(n / \nu)|^2 \ll \sum_{1 \le n \le \nu H^{\epsilon / 4}} 1 + \sum_{n > \nu H^{\epsilon / 4}} \frac{H^{\epsilon / 2} \nu^2}{n^2} \ll_\epsilon \nu H^{\epsilon / 4}
\]
for any $\nu > 0$. Hence,
\[
2 H^2 \sum_{b > X^{1/4}} \frac{\mu^2(b)}{b^3} \sum_{n \ge 1} \Big| w\Bigl( \frac{n H}{b^{3/2}} \Bigr) \Big|^2 \ll H^{1 + \epsilon / 4} \sum_{b > X^{1/4}} \frac{1}{b^{3/2}} \ll_\epsilon \frac{H^{1 + \epsilon/4}}{X^{1/8}}.
\]
Thus, the left hand side of \eqref{smoothW} is
\begin{equation} \label{smoothW2}
2 H^2 \sum_{b \ge 1} \frac{\mu^2(b)}{b^3} \sum_{n \ge 1} \Big| w\Bigl( \frac{n H}{b^{3/2}} \Bigr) \Big|^2 + O_\epsilon \Bigl(\frac{H^{1 + \epsilon/4}}{X^{1/8}}\Bigr).
\end{equation}

Secondly, we use contour integration to simplify \eqref{smoothW2}. Define $g(x) := |w(e^x)|^2 e^x$, which is smooth and decays exponentially as $|x| \rightarrow \infty$ by \eqref{deriv}. Its Fourier transform is given by
\[
\hat{g}(\xi) = \int_{-\infty}^{\infty} |w(e^x)|^2 e^x e^{-2 \pi i x \xi} dx = \int_{0}^{\infty} |w(y)|^2 y^{-2 \pi i \xi} dy,
\]
and one can show that (i) $\hat{g}(\xi)$ is entire and (ii) $\hat{g}(\xi) = O(H^{2\epsilon} / (|\xi| + 1)^3)$ uniformly for $|\Im(\xi)| < 1 / (2 \pi)$ by \eqref{deriv} and standard partial summation argument. Then Fourier inversion gives, for $r > 0$,
\[
|w(r)|^2 = r^{-1} \frac{1}{2 \pi i} \int_{(c)} r^s \hat{g} \Bigl( \frac{s}{2 \pi i} \Bigr) ds,
\]
where the integral is over $\Re(s) = c$ with $-1 < c < 1$. Hence, taking $c = - 1/4$,
\[
2 H^2 \sum_{b \ge 1} \frac{\mu^2(b)}{b^3} \sum_{n \ge 1} \Big| w\Bigl( \frac{n H}{b^{3/2}} \Bigr) \Big|^2 = \frac{H}{i \pi} \sum_{b \ge 1} \frac{\mu^2(b)}{b^{3/2}} \sum_{n \ge 1} \frac{1}{n} \int_{(-1/4)} H^s n^s b^{-3s/2} \hat{g}\Bigl( \frac{s}{2 \pi i} \Bigr) ds.
\]
Due to absolute convergence of the series, one can take the sums inside the integral and the above becomes
\[
\frac{H}{i \pi} \int_{(-1/4)}  H^s  \zeta(1 - s) \frac{\zeta(\frac{3}{2}(1 + s))}{\zeta(3(1+s))}  \hat{g}\Bigl( \frac{s}{2 \pi i} \Bigr) ds
\]
by \eqref{zeta}. Because of the bounds on $\hat{g}$ and $\zeta(1/2 + it) \ll t^{1/6} \log^2 t$, we can push the contour integral above to the left to an integral over $\Re(s) = -2/3$. Picking out a residue at $s = -1/3$ from the singularity of $\zeta(\frac{3}{2}(1 + s))$, the above simplifies to
\[
\frac{4 \zeta(4/3)}{3 \zeta(2)} \hat{g} \Bigl( - \frac{1}{6 \pi i} \Bigr) H^{2/3} + O(H^{1/3 + 2\epsilon}) =  \frac{4 \zeta(4/3)}{3 \zeta(2)} \int_{0}^{\infty} |w(y)|^2 y^{1/3} dy \cdot H^{2/3} + O(H^{1/3 + 2 \epsilon}).
\]
This and \eqref{smoothW2} finish the proof of the claim.

\bigskip

To deal with \eqref{diagonal-main}, we introduce the function
\[
w(y) := S(y) h\Bigl( \frac{y}{H^{\epsilon/4}} \Bigr)
\]
where $h$ is a smooth bump function such that $h(x) = 1$ for $|x| \le 1$ and $h(x) = 0$ for $|x| \ge 2$ (for example, one can take $h(x) = \sigma_{10, 2, 0.1}^{+}(x + 3)$). Then one can check that $w(y)$ satisfies the hypothesis \eqref{deriv}. For such $w$,
\begin{align} \label{sum-error}
2H^2 \sum_{b \le X^{1/4}} & \frac{\mu^2(b)}{b^3} \sum_{n \ge 1} \biggl( S\Bigl( \frac{n H}{b^{3/2}} \Bigr)^2 - w \Bigl( \frac{n H}{b^{3/2}} \Bigr)^2 \biggr) \ll H^2 \mathop{\sum_{b \le X^{1/4}, n \ge 1}}_{n H^{1 - \epsilon/4} \ge b^{3/2}} \frac{1}{b^3} \frac{1}{(n H / b^{3/2})^2} \nonumber \\
&\ll \sum_{b \le H^{2/3 - \epsilon/6}} \sum_{n \ge 1} \frac{1}{n^2} + \sum_{b > H^{2/3 - \epsilon/6}} \sum_{n \ge b^{3/2} / H^{1 - \epsilon/4}} \frac{1}{n^2} \ll_\epsilon H^{2/3 - \epsilon/6}.
\end{align}
On the other hand,
\begin{equation} \label{integral-error}
\int_{0}^{\infty} |S(y)|^2 y^{1/3} dy - \int_{0}^{\infty} |w(y)|^2 y^{1/3} dy \ll \int_{H^{\epsilon/4}}^{\infty} y^{-5/3} dy \ll_\epsilon H^{-\epsilon / 6}.
\end{equation}
Combining \eqref{smoothW}, \eqref{sum-error} and \eqref{integral-error}, we have
\[
2 H^2 \sum_{b \le X^{1/4}} \frac{\mu^2(b)}{b^3} \sum_{n \ge 1} \Big| S\Bigl( \frac{n H}{b^{3/2}} \Bigr) \Big|^2 = H^{2/3} \frac{4 \zeta(4/3)}{3 \zeta(2)} \int_{0}^{\infty} |S(y)|^2 y^{1/3} dy + O_\epsilon \Bigl(\frac{H^{1 + \epsilon/4}}{X^{1/8}} + H^{2/3 - \epsilon/6}\Bigr)
\]
which gives Proposition \ref{prop2}.


\section{Proof of Theorem \ref{thm1}}

From Propositions \ref{prop1} and \ref{prop2}, we have
\[
\frac{1}{X} \int_{X}^{2X} \Big| Q(x+y) - Q(x) - \frac{\zeta(3/2)}{\zeta(3)} H \Big|^2 dx \ll_\epsilon H^{4/3} + \frac{H^{5/3 + \epsilon}}{X^{1/16}} + \frac{H^{2+\epsilon}}{X^{1/8}} \ll H^{16/11 + \epsilon}
\]
for $X^{3 \epsilon} \le H \le X^{11 / 48}$. For $X^{11 / 48} < H \le X^{1/2 - \epsilon}$, we simply apply \eqref{fullcount} and get
\begin{align*}
\frac{1}{X} \int_{X}^{2X} \Big| Q(x+y) - Q(x) - \frac{\zeta(3/2)}{\zeta(3)} H \Big|^2 dx \ll& \frac{1}{X} \int_{X}^{2X} \Big|\frac{\zeta(2/3)}{\zeta(2)} (\sqrt[3]{x+y} - \sqrt[3]{x}) + x^{1/6} \Big|^2 dx \nonumber \\
\ll& \frac{H^2}{X^{1/3}} + X^{1/3} \ll H^{16 / 11}.
\end{align*}
Combining the above, we have Theorem \ref{thm1}.

\section{A sketch of proof of Theorem \ref{thm2}}

Similar to our proof of Proposition \ref{prop1}, one starts by splitting the integral \eqref{z1} into
\[
\frac{1}{X} \int_{X}^{2X} \bigg| \mathop{\sum_{x < a^2 b^3 \le x + y}}_{b \le H^{2/3 + \epsilon}} \mu^2(b) - H \sum_{b \le H^{2/3 + \epsilon}} \frac{\mu^2(b)}{b^{3/2}} + \mathop{\sum_{x < a^2 b^3 \le x + y}}_{b > H^{2/3 + \epsilon}} \mu^2(b) - H \sum_{b > H^{2/3 - \epsilon}} \frac{\mu^2(b)}{b^{3/2}} \bigg|^2 dx.
\]
It can be rewritten as $I_1 + I_2 + O( I_1^{1/2} I_2^{1/2} )$ where
\[
I_1 = \frac{1}{X} \int_{X}^{2X} \bigg| \mathop{\sum_{x < a^2 b^3 \le x + y}}_{b \le H^{2/3 + \epsilon}} \mu^2(b) - H \sum_{b \le H^{2/3 + \epsilon}} \frac{\mu^2(b)}{b^{3/2}} \bigg|^2 dx
\]
and
\[
I_2 = \frac{1}{X} \int_{X}^{2X} \bigg| \mathop{\sum_{x < a^2 b^3 \le x + y}}_{b > H^{2/3 + \epsilon}} \mu^2(b) - H \sum_{b > H^{2/3 + \epsilon}} \frac{\mu^2(b)}{b^{3/2}} \bigg|^2 dx.
\]
Then, one can imitate the proof of Proposition 2 in \cite{GMRR} by splitting $I_2$ into dyadic ranges and trying to show that
\begin{equation} \label{newmean}
\frac{1}{X} \int_{X}^{2X} \bigg| \mathop{\sum_{x < a^2 b^3 \le x + y}}_{B < b \le 2B} \mu^2(b) - H \sum_{B < b \le 2B} \frac{\mu^2(b)}{b^{3/2}} \bigg|^2 dx \ll_\epsilon H^{2/3 - \epsilon/4}
\end{equation}
for $B \in [H^{2/3 + \epsilon}, (2X)^{1/3}]$. Following \cite{GMRR}, one can use contour integral of Dirichlet series and Plancherel's identity to bound the left hand side of \eqref{newmean} by
\[
\ll \sqrt{X} \int_{-\infty}^{\infty} \min \Bigl( \Bigl( \frac{H}{\sqrt{X}} \Bigr)^2, \frac{1}{|t|^2} \Bigr) \, |\zeta(1/2 + 2 i t) M( 3/4 + 3 i t ) |^2 dt \; \; \text{ with } \; \; M(s) := \sum_{B < b \le 2 B} \frac{\mu^2(b)}{b^s}.
\]
One can show that the integral with $|t| \ge X$ contributes $O(1)$. For the remaining range of $t$, one can split the range and show that it is
\[
\ll H \Bigl( \sup_{\sqrt{X} / H \le T \le X} \frac{1}{T} \int_{|t| \le T} |\zeta(1/2 + 2 i t) M( 3/4 + 3 i t ) |^2 dt \Bigr) + O(1).
\]
Applying Lindel\"{o}f Hypothesis and standard mean-value theorem on Dirichlet polynomial, the above can be shown to be
\[
\ll \frac{H T^\delta}{T} (T + B) \cdot \frac{1}{\sqrt{B}} \ll \frac{H T^\delta}{\sqrt{B}} + \frac{H T^\delta \sqrt{B}}{T}
\]
for some $\delta > 0$. As $\sqrt{X} / H \le T \le X$ and $B \gg H^{2/3 + \epsilon}$, the above is $O_\epsilon(H^{2/3 - \epsilon/4})$ when $H \ll X^{1/4 - \epsilon}$.

\bigskip

It remains to deal with $I_1$. This can be done similarly as our estimation of $J_1$ using smooth weight. The only difference is in some of the details in our counting argument for the non-diagonal terms. We have $B_1, B_2 \le H^{2/3 + \epsilon}$. Instead of \eqref{bound3}, the non-diagonal contribution is
\[
\ll \sum_{b} \sum_{n} \sum_{b_1'} \sum_{b_2'} \frac{1}{N_1 N_2} \cdot \frac{B_1^{3/2} N_2}{n^2 X^{1/2 - \epsilon}} + \sum_{b} \sum_{n} \sum_{b_1'} \sum_{b_2'} \frac{1}{N_1 N_2} \ll \frac{B_1^{5/2} B_2}{X^{1/2 - \epsilon}} + \sum_n \frac{B_1 B_2}{N_1 N_2}
\]
which is $O_\epsilon(H^{2/3 - \epsilon})$ as long as $H \le X^{3/14 - 5 \epsilon}$ and $n \le N_1 N_2 / H^{2/3 + 3 \epsilon}$. It remains to deal with $n > N_1 N_2 / H^{2/3 + 3 \epsilon}$. One can show that $b \le (B_1 B_2)^{1/2 - \delta}$ for some $\delta > 0$ when $H \ll X^{3/14 - 2 \delta - 5 \epsilon}$. Then, instead of \eqref{bound4}, the non-diagonal contribution is
\[
\sum_{b} \sum_{n} \sum_{n_1'} \sum_{n_2'} \frac{1}{N_1 N_2} \cdot \frac{B_1 B_2^{5/2}}{b^2 N_2 X^{1/2 - \epsilon}} + \sum_{b \le (B_1 B_2)^{1/2 - \delta}} \sum_{n} \sum_{n_1'} \sum_{n_2'} \frac{1}{N_1 N_2} \ll \frac{B_1 B_2^{5/2}}{X^{1/2 - \epsilon}} + (B_1 B_2)^{1/2 - \delta}
\]
which is $O_\epsilon(H^{2/3 - \epsilon})$ as long as $H \ll X^{3/14 - 2 \delta - 5 \epsilon}$ and $\delta \ge \epsilon$. So, we have an acceptable error.

\bigskip

Finally, the main term of $I_1$ can be studied in the same way as the proof of Proposition \ref{prop2} with an acceptable error $O_\epsilon(H^{2/3 - \epsilon/6})$.

\bigskip

{\bf Acknowledgment:} The author would like to thank the anonymous referee for pointing out several inaccuracies and helpful suggestions which led to modified arguments and improved results. He also wants to thank Ofir Gorodetsky for pointing out Edva Roditty-Gershon's related work on squarefull polynomials.

\bibliographystyle{amsplain}

Mathematics Department \\
Kennesaw State University \\
Marietta, GA 30060 \\
tchan4@kennesaw.edu

\end{document}